\documentclass[12pt,reqno]{amsart}
\newtheorem{thm}{Theorem}[section]
\newtheorem{defi}{Definition}[section]

\newtheorem{cor}{Corollary}[section]
\newtheorem{pr}{Proposition}[section]
\theoremstyle{definition}
\newtheorem{rem}{Remark}[section]

\newcommand{\be}{\begin{equation}}
\newcommand{\ee}{\end{equation}}
\newcommand{\bea}{\begin{eqnarray}}
\newcommand{\eea}{\end{eqnarray}}
\newcommand{\beb}{\begin{eqnarray*}}
\newcommand{\eeb}{\end{eqnarray*}}
\usepackage{amssymb,amsfonts,amsthm,setspace,indentfirst,mathrsfs}
\usepackage{minibox}
\usepackage{enumitem}
\numberwithin{equation}{section}
\begin{document}
%
\title{  Curvature properties of Bardeen black hole spacetime }
\author[A. A. Shaikh., S. K. Hui. and M. Sarkar]{Absos Ali Shaikh, Shyamal kumar hui, Mousumi Sarkar}
\date{\today}
\address{\noindent\newline$^1$ Department of Mathematics,
	\newline University of Burdwan, 
	\newline Golapbag, Burdwan-713104,
	\newline West Bengal, India} 
\email{aask2003@yahoo.co.in, aashaikh@math.buruniv.ac.in}
\email{skhui@math.buruniv.ac.in}
\email{mousrkr13@gmail.com}

\dedicatory{}
\begin{abstract}
The Bardeen solution corresponding to Einstein field equations with a cosmological constant is a regular black hole. The main  goal of this manuscript is to investigate the geometric structures in terms of curvature conditions admitted by this spacetime. It is found that this spacetime is pseudosymmetric and possess several kinds of pseudosymmetries. Also, it is a manifold of pseudosymmetry Weyl curvature and the difference tensor $C\cdot R - R\cdot C$ linearly depends on the tensors $Q(g,C)$ and $Q(S,C)$. It is interesting to note that such a spacetime is weakly generalized recurrent manifold and satisfies special recurrent like structure. Further, it is an Einstein manifold of level $2$  and Roter type. The energy momentum  tensor of this spacetime is pseudosymmetric and finally a worthy comparison between the geometric properties of Bardeen spacetime and Reissner-Nordstr\"om spacetime is given.
\end{abstract}
%
\maketitle
%

\section{\bf Introduction}\label{intro}
%
\indent The geometry of a space is described by the curvature, which plays a crucial role in differential geometry as this symmetry of the space is determined by the restriction on the curvature tensor $R$. There are several classes of manifolds with specific geometric structures such as the locally symmetric manifolds by Cartan \cite{Cart26}  defined as $\nabla R= 0$, semisymmetric manifolds again by Cartan \cite{Cart46,Szab82,Szab84,Szab85}  defined as $R \cdot R=0$,  pseudosymmetric manifolds by Adam\'{o}w and Deszcz \cite{AD83}. Many authors have studied locally symmetric manifolds and introduced several notions of manifolds, such as recurrent manifolds by Ruse \cite{Ruse46, Ruse49a, Ruse49b, Walk50},  several classes of generalized recurrent manifolds by Shaikh and his coauthors \cite{SP10, SR10,SRK18,SRK17},  weakly symmetric manifolds by Tam\'assy and Binh \cite{TB89, TB93}, pseudosymmetric manifolds by Chaki \cite{Chak87, Chak88}, curvature $2$-forms of recurrent manifolds by Besse \cite{Bess87, LR89, MS12a,MS13a,MS14} etc. We mentioned that the notion of pseudosymmetry in the sence of Deszcz is important in pseudo-Riemannian geometry as well as in general relativity because several spacetimes are models of pseudosymmetric manifolds (see e.g. \cite{ ADEHM14, DK99, Kowa06, SAA18, SAAC20,  SDKC19}). It is mention that  pseudosymmetry in the sence of Deszcz and pseudosymmetry in the sence of Chaki are different concepts. The study of geometric structure of a certain spacetime gives us the idea about its geometry and provides knowledge about its physical nature.\\
\indent It is known that spacetime singularities are a reflection of the incompletenes of general relativity.  In $1968$,  Bardeen \cite{Bardeen} obtained a solution of Einstein field equation in spherical symmetry, which describes a singularity free black hole spacetime.  In fact, it is the first regular black hole model in general relativity. In spherical coordinates system $(t,r,\theta, \phi)$, the Bardeen spacetime  metric is as follows   
\bea\label{BM}
ds^2=-\left( 1-\frac{2M\rho^2}{(e^2+\rho^2)^{3/2}}\right)dt^2  +\left(1-\frac{2M\rho^2}{(e^2+\rho^2)^{3/2}}\right)^{-1} dr^2 + r^2\left( d\theta^2 + \sin^2 \theta d\phi^2\right),  
\eea
where, $M$ and $e$ respectively denotes the mass and  magnetic charge of the nonlinear self gravities monopole \cite{Bardeen}.
If  $e^2\le \frac{16}{27}M^2$, then the Bardeen model represents a black hole and a singularity-free structure \cite{AE}. When $e^2 = \frac{16}{27}M^2$, the horizons shrink into a single one that corresponds to an extreme black hole such as the extreme Reissner-Nordstr\"{o}m solution. In \cite{Borde}, Borde has  studied Bardeen spacetime and proved that topology change is necessary for the existence of regular black holes satisfying the weak energy condition. In \cite{AE} Ay\'{o}n-Beato and Garc\'{i}a  reinterpreted this black hole as a gravitational field of a nonlinear magnetic monopole. In the context of black holes we also mention the works of Capozziello and his coauthors \cite{BCL03, BCIS01, CFS04} regarding its gravitational lensing, cosmological structure, types as part of entangled systems etc. 
\\
\indent  However, the curvature properties of Bardeen spacetime are yet to known. Thus the present study is devoted to deduce the geometric structures of Bardeen black hole metric in terms of curvature restrictions.

\indent Deriving the  components of various curvature tensors, we  investigate the curvature properties admitted by the Bardeen black hole metric \eqref{BM}. It is found that Bardeen black hole metric is not semisymmetric but  pseudosymmetric manifold and satifies various pseudosymmetric type curvature conditions. The difference tensor $C\cdot R -R\cdot C$ is linearly dependent with the tensors $Q(g,C)$ and $Q(S,C)$. It is also a weakly generalized recurrent manifold and satisfies special recurrent like structure. Further, it is Roter type and Einstein manifold of level $2$. It is interesting to note that this spacetime has pseudosymmetric energy momentum tensor. \\   
\indent In this context, we mention that recently the   curvature properties of interior black hole metric \cite{SDHK} were investigated and it is shown that such a interior black hole metric admits several geometric structures such as pseudosymmetry, Ricci pseudosymmetry and several kinds of pseudosymmetric type curvature conditions.\\ 
\indent This paper is composed in four sections : in section $2$ we review  various curvature tensors and  geometric structures as priliminaries.  Section $3$  determines the  curvature restricted geometric structures admitted by the Bardeen black hole metric \eqref{BM}. Section $4$ deals with some geometric structures of energy momentum tensor. In section $5$, we make a comparison between the geometric properties of Bardeen black hole spacetime and Reissner-Nordstr\"om spacetime.

%

%
\section{\bf Preliminaries}
\indent Let $M$, a smooth and connected manifold of dimension $\ge 3$, be furnished with the Levi-Civita connection $\nabla$ and also with a semi-Riemannian metric $g$. Let $\kappa$ be the scalar curvature, $S$ the Ricci curvature and $R$ the Riemann curvature of $M$.\\
 Now, for a symmetric $(0,2)$  tensors  $\lambda$ and $\tau$, the Kulkarni-Nomizu product $\lambda \wedge \tau$ \cite{DGHS11, Glog02, SKA16} is defined as 
\bea
(\tau \wedge  \lambda )(\zeta_1,\zeta_2,X,Y)&=& \tau  (\zeta_1,Y)\lambda (\zeta_2,X) - \tau (\zeta_1,X)\lambda(\zeta_2,Y) \nonumber\\
&+& \tau (\zeta_2,X)\lambda(\zeta_1,Y)  - \tau (\zeta_2,Y)\lambda(\zeta_1,X)\nonumber
\eea
$\zeta_1, \zeta_2, X, Y$ being some smooth vector fields and throughout this paper we assume  $\zeta_1,\zeta_2, X,Y \in  \chi({M})$, $\chi(M)$ being the set of all smooth vector fields of $M$.\\
\indent Now, we consider the following endomorphisms on $M$ (\cite{DHJKS14,SC21,SK18,SK19}) :
\bea
(\zeta_1\wedge_{\tau} \zeta_2)X &=& \tau(\zeta_2,X)\zeta_1-\tau(\zeta_1,X)\zeta_2 ,\nonumber \\
\mu_{\mathscr{R}}(\zeta_1,\zeta_2)X &=& \left( [\nabla_{\zeta_1}, \nabla_{\zeta_2}] - \nabla_{[\zeta_1,\zeta_2]}\right) X, \nonumber \\
\mu_{\mathscr{W}}(\zeta_1,\zeta_2)X &=& \mu_{\mathscr{R}}(\zeta_1,\zeta_2)X-\frac{\kappa}{n(n-1)}\left( \zeta_1\wedge_g\zeta_2\right)X ,\nonumber \\
\mu_{\mathscr{K}}(\zeta_1,\zeta_2)X &=& \mu_{\mathscr{R}}(\zeta_1,\zeta_2)X -\frac{1}{(n-2)}\left( \zeta_1\wedge_S\mathscr{J}\zeta_2 +\mathscr{J}\zeta_1\wedge_S\zeta_2\right)X,  \nonumber \\
\mu_{\mathscr{C}}(\zeta_1,\zeta_2)X &=& \mu_{\mathscr{R}}(\zeta_1,\zeta_2)X +\frac{\kappa}{(n-1)(n-2)}(\zeta_1\wedge_g\zeta_2)X \nonumber\\ &-& \frac{1}{(n-2)}\left(\zeta_1\wedge_S\mathscr{J}\zeta_2 +\mathscr{J}\zeta_1\wedge_S\zeta_2\right)X , \nonumber \\
\mu_{\mathscr{P}}(\zeta_1,\zeta_2)X &=& \mu_{\mathscr{R}}(\zeta_1,\zeta_2)X-\frac{1}{(n-1)}\left( \zeta_1\wedge_S\zeta_2\right)X  ,\nonumber 
\eea
\indent where $\mathscr{J}$ is the Ricci operator and is defined by $g(\zeta_1,\mathscr{J}\zeta_2)=S(\zeta_1,\zeta_2)$.\\
\indent Now, for an endomorphism $\mathscr{D}(\zeta_1,\zeta_2)$, the corresponding  $(0,4)$-tensor $D$ is given by
\bea
D(\zeta_1,\zeta_2,\zeta_3,\zeta_4)=g(\mathscr{D}(\zeta_1,\zeta_2)\zeta_3,\zeta_4).
\eea
In above relation replacing $\mathscr{D}$ by $\mu_{\mathscr{C}}$ one can find the $(0,4)$ conformal curvature $C$. Again, replacing it by $\mu_{\mathscr{R}}$, $\mu_{\mathscr{P}}$, $\mu_{\mathscr{W}}$ and $\mu_{\mathscr{K}}$ one can also look for $(0,4)$ Riemann curvature $R$, projective curvature $P$, concircular curvature $W$ and conharmonic curvature $K$ respectively.

\indent Let $\eta$ be a $(0,l)$-tensor field on  $M$, $l\geq 1$. Now we can define the $(0,l+2)$-tensor $D \cdot \eta$ \cite{DGHS98, DH03, SK14} by 
\beb
(D\cdot \eta)(\zeta_1,\zeta_2,\cdots,\zeta_l;X,Y)&=&(\mathscr{D}(X,Y)\cdot \eta)(\zeta_1,\zeta_2,\cdots,\zeta_l)\\
&&\hspace{-1in} =-\eta(\mathscr{D}(X,Y)\zeta_1,\zeta_2,\cdots,\zeta_l)- \cdots -\eta(\zeta_1,\zeta_2,\cdots,\mathscr{D}(X,Y)\zeta_l),
\eeb
\index In addition, if $\lambda$ is a symmetric $(0,2)$-tensor field then we define the $(0,l+2)$-tensor $Q(\lambda, \eta)$, called Tachibana tensor \cite{DGPSS11,SDHJK15,Tach74}, by 
\beb
Q(\lambda,\eta)(\zeta_1,\zeta_2,\cdots,\zeta_l;X,Y)&=&((X\wedge_{\lambda}Y)\cdot \eta)(\zeta_1,\zeta_2,\cdots,\zeta_l)\\
&&\hspace{-2in} =-\eta((X\wedge_{\lambda}Y)\zeta_1,\zeta_2,\cdots,\zeta_l)- \cdots -\eta(\zeta_1,\zeta_2,\cdots,(X\wedge_{\lambda}Y)\zeta_l)\\
&&\hspace{-2in} = \lambda(\zeta_1, X) \eta(Y,\zeta_2,\cdots,\zeta_l) + \cdots + \lambda(\zeta_l,X) \eta(\zeta_1,\zeta_2,\cdots,Y)\\
&&\hspace{-2in} - \lambda(\zeta_1,Y) \eta(X,\zeta_2,\cdots,\zeta_l) - \cdots - \lambda(\zeta_l,Y) \eta(\zeta_1,\zeta_2,\cdots,X).
\eeb
\index The local components $(D \cdot \eta)_{i_1i_2 \cdots i_l\alpha \beta}$ and $Q(\lambda, \eta)_{i_1i_2 \cdots i_l\alpha \beta}$ of the tensors $D \cdot \eta $ and the Tachibana tensor  $Q(\lambda,\eta)$ are written as : 
\beb
(D\cdot \eta)_{i_1i_2...i_l\alpha \beta} &=& -g^{uv}[D_{\alpha \beta i_1v}\eta_{ui_2...i_l} +\cdots + D_{\alpha \beta i_lv}\eta_{i_1i_2...u}],\\
Q(\lambda,\eta)_{i_1i_2...i_l\alpha \beta} &=& \lambda_{i_1 \beta }\eta_{\alpha i_2...i_l} + \cdots + \lambda_{i_l \beta }\eta_{i_1i_2...\alpha} \\ &-& \lambda_{i_1 \alpha}\eta_{\beta i_2...i_l} - \cdots - \lambda_{i_l \alpha}\eta_{i_1i_2...\beta}.
\eeb
\begin{defi}
	The linear dependency of the tensor $D\cdot \eta$ with $Q(\lambda, \eta)$ defines a $\eta$-pseudosymmetric manifold \cite{AD83,  Desz92, Desz93, DGHZ15, DGHZ16, SK14, SKppsnw} due to $D$ i.e., on this manifold we have $D\cdot \eta=\varrho _\eta Q(\lambda, \eta)$ where $\varrho_\eta$ is a smooth function on $\left\{ x \in M : Q(\lambda,\eta)_x \ne 0 \right\}$. If $D\cdot \eta=0$ holds then the manifold $M$ is a $\eta$-semisymmetric manifold due to $D$ \cite{Cart46,SAAC20N,Szab82, Szab84, Szab85}.
\end{defi}
For $D=R$, $\eta=R$ and $\lambda=g$ (or $S$) the manifold is simply called as pseudosymmetric manifold or Deszcz's pseudosymmetric manifold (or Ricci generalized pseudosymmetric manifold). Again a semisymmetric manifold is obtained for $D=R$, $\eta=R$ i.e., $R\cdot R=0$.  Several kinds of pseudosymmetric and semisymmetric type curvature conditions are obtained for other curvature tensors. \\
\indent If $S=\alpha g$ holds on $M$ i.e., the Ricci tensor and the metric tensor are linearly dependent then it is called an Einstein manifold and in this case $\alpha=\frac{\kappa}{n}$. In a quasi-Einstein manifold the Ricci tensor has the form $S=\beta A\otimes A + \alpha g$ (A being some $1$-form). This notion is important in general relativity as a connected Lorentzian manifold of dimension $4$ is a spacetime of perfect fluid if it is quasi-Einstein and vise versa. In a quasi-Einstein manifold rank$(S-\alpha g)=1$ and if rank$(S-\alpha g)=2\ \text{or}\ 3$ then the manifold is called as $2$ or $3$ quasi-Einstein manifold \cite{S09, SK19} respectively. It may be mentioned that Kaigorodov spacetime \cite{SDKC19} is Einstein, Robertson Walker spacetime \cite{ARS95, Neill83, SKMHH03} are quasi Einstein and Kantowski-Sachs spacetime \cite{SC21} is $2$-quasi Einstein. Again, a Ricci simple manifold is a special case of quasi Einstein manifold for $\alpha=0$ and Morris-Thorne spacetime \cite{ECS22} is a Ricci simple manifold.  Another generalization of Einstein manifold is given  below. \\
\begin{defi}$($\cite{Bess87, SK14, SK19}$)$
		A manifold satisfying the relation
		\be
		\beta_1S^4 + \beta_2S^3 + \beta_3S^2 + \beta_4S + \beta_5g=0, \; (\beta_1 \neq 0)\nonumber
		\ee
		is called an Einstein manifold of level $4$, 
		where $\beta_i$ are smooth functions on $M$. If $\beta_1=0$ but $\beta_2\ne 0$ (resp., $\beta_1=\beta_2=0$ but $\beta_3 \ne 0$), then it turns into an Einstein manifold of level $3$ (resp. level $2$).
\end{defi}
\begin{defi}
	A manifold is called a generalized Roter type manifold \cite{Desz03, DGJPZ13, DGJZ-2016, DGP-TV-2015, SK16,SK19} if its  Riemann curvature tensor fulfills  $$R=s_{22}(g\wedge g)+ (s_{11}S+s_{12}g)\wedge S+ (s_{00}S^2+s_{01}S+s_{02}g)\wedge S^2$$  $s_{ij}$ being some smooth functions on $M$. For $s_{00}=s_{01}=s_{02}=0$, $R$ is linearly dependent on $g\wedge g$, $S\wedge g$ and $S\wedge S$ and in this case $M$ is called a Roter type manifold \cite{Desz03, DG02, DGP-TV-2011, DPSch-2013, Glog-2007}. 
\end{defi}
\indent It is worthy to note that Melvin magnetic spacetime \cite{SAAC20} and Robinson-Trautman spacetime \cite{SAA18} are Roter type and $Ein(2)$ manifolds, whereas Lifshitz spacetime \cite{SSC19} is $Ein(3)$ manifold and generalized Roter type. Recently, Lema\^itre-Tolman-Bondi spacetime \cite{SAACD22} is found to possess generalized Roter type and $Ein(3)$ condition.
\begin{defi}
	A super generalized recurrent manifold $M$ is defined by curvature condition \cite{SRK17}
	$$\nabla R=\Pi\otimes R + A\otimes (g\wedge g) + \bar{A} (S\wedge g) + \bar{\bar{A}}(S\wedge S)$$
	where $A, \bar{A}, \bar{\bar{A}}$ are some $1$-forms on $M$. A weakly generalized recurrent manifold \cite{SAR13,SR11} (resp., hyper generalized recurrent manifold \cite{SP10,SRK18}) is obtained for $A=\bar{A}=0$ (resp., $A=\bar{\bar{A}}=0$).
\end{defi}
\begin{defi}
	Tam\'assy and Binh \cite{TB89, TB93}  defined the notion of weak symmetry by 
	\beb
	(\nabla_{X} R)(\zeta_1,\zeta_2,\zeta_3,\zeta_4)&=& \Pi(X)\otimes R(\zeta_1,\zeta_2,\zeta_3,\zeta_4)+ \bar{B}(\zeta_4)\otimes R(\zeta_1,\zeta_2,\zeta_3,X)\\ &+& B(\zeta_3)\otimes R(\zeta_1,\zeta_2,X,\zeta_4)+ \bar{A}(\zeta_2)\otimes R(\zeta_1,X,\zeta_3,\zeta_4)\\ &+& A(\zeta_1)\otimes R(X,\zeta_2,\zeta_3,\zeta_4)
	\eeb
	where $A, \bar{A}, B, \bar{B}$ are some $1$-forms on $M$. For $\Pi=\frac{A}{2}=\frac{\bar A}{2}=\frac{B}{2}=\frac{\bar B}{2}$ it is a Chaki pseudosymmetric manifold \cite{Chak87, Chak88}.
\end{defi}
\begin{defi}\label{defi2.8} \cite{LR89, MS12a, MS13a, MS14}
	Let $D$ be a  $(0,4)$-type tensor and $\lambda$ be a $(0,2)$-type symmetric tensor on $M$. Then $\Omega^m_{(D)l}$ \cite{SKP03}, the curvature $2$-forms, are recurrent if 
	\beb
	\mathop{\mathcal{S}}_{\zeta_1,\zeta_2,\zeta_3}(\nabla_{\zeta_1}D)(\zeta_2,\zeta_3,X,Y)=\mathop{\mathcal{S}}_{\zeta_1,\zeta_2,\zeta_3}A (\zeta_1)\otimes D(\zeta_2,\zeta_3,X,Y)
	\eeb
	  $\mathcal S$ being the cyclic sum over $\zeta_1$, $\zeta_2$ and $\zeta_3$, holds on $M$ and the $1$-forms $\iota_{(\lambda)l}$ \cite{SKP03} are recurrent if 
	$$(\nabla_{\zeta_1}\lambda)(\zeta_2,X)-(\nabla_{\zeta_2}\lambda)(\zeta_1,X)=A(\zeta_1)\otimes \lambda(\zeta_2,X)-\bar{A}(\zeta_2)\otimes \lambda(\zeta_1,X)$$ holds on $M$ and $A$, $\bar A$ are smooth 1-forms.
\end{defi}
\begin{defi}\label{2.6}
The Ricci tensor of $M$ is of Codazzi type \cite{F81, S81} (resp., cyclic parallel \cite{ Gray78, SB08, SJ06, SS06, SS07}) if 
	$$(\nabla_{\zeta_1}S)(\zeta_3,\zeta_2)=(\nabla_{\zeta_2}S)(\zeta_3,\zeta_1)$$
	$$(resp., \mathop{\mathcal{S}}_{\zeta_1,\zeta_2,\zeta_3}(\nabla _{\zeta_1} S)(\zeta_2,\zeta_3)=0)$$
	holds on $M$,
\end{defi}
\indent We mention that the Ricci tensor of the  $(t-z)$-type plane wave metric is of Codazzi type \cite{EC21} and the Ricci tensor is cyclic parallel in G\"odel spacetime \cite{DHJKS14}. 
\begin{defi}\label{defi2.7}$($\cite{DD91, DGJPZ13, MM12a, MM12b, MM13}$)$
	The Ricci tensor $S$ of $M$ is called  Riemann compatible if the relation   $$\mathop{\mathcal{S}}_{\zeta_1,\zeta_2,\zeta_3}R(\mathscr J  \zeta_1,X,\zeta_2,\zeta_3)=0$$ holds.
\end{defi}
\indent This notion of compatibility can be extended to the curvatures $C$, $P$, $W$ and $K$ to define the corresponding curvatures compatibility.
\begin{defi}$($\cite{P95, Venz85}$)$
	Let $D$ be a symmetric $(0,4)$-type tensor of $M$. If the $1$-forms $\Theta$ satisfying the relation 
	\beb
	\mathop{\mathcal{S}}_{\zeta_1,\zeta_2,\zeta_3}\Theta(\zeta_1)\otimes D(\zeta_2,\zeta_3,X,Y)=0
	\eeb
	generates a $k$-dimensional vector space, $k \ge 1$, then $M$ is a $D$-space by Venzi.
	\end{defi}
%
\section{\bf Geometric properties admitted by Bardeen black hole metric}
%
The components of metric \eqref{BM} are

$$\begin{array}{c}g_{11}=-\left(1-\frac{2M\rho^2}{\rho^3_1}\right),\; g_{22}= \left((1-\frac{2M\rho^2}{\rho^3_1})^{-1}\right),\;  \\ g_{33}=\rho^2, \;g_{44}=\rho^2 \sin^2\theta, \; g_{ij}=0, \;\mbox{otherwise}.
\end{array}$$

\index Now, the components of various curvature tensors of Bardeen black hole metric \eqref{BM} are calculated in a straight forward manner. \\
\index The non-zero components of second kind Christoffel symbols $\Gamma^h_{ij}$ of $g$ are given by : 
\begin{eqnarray}
\begin{cases}
\Gamma^2_{11}=-\frac{M\rho (\rho^2_1-3e^2) (2M\rho^2- \rho^3_1)}{\rho^8}; \ \
\Gamma^1_{12} = \frac{M\rho (\rho^2_1-3e^2)}{\rho^2_1(-2M\rho^2+{\rho^3_1})}=-\Gamma^2_{22} ; \\
\Gamma^3_{23} = \frac{1}{\rho}=\Gamma^4_{24};  \\
\Gamma^2_{33}=-\rho+\frac{2M\rho^3}{\rho^3_1}; \\
\Gamma^4_{34}=\cot\theta; \\
\Gamma^2_{44}=\rho[-1+\frac{2M\rho^2}{\rho^3_1}]\sin^2\theta; \\
 \Gamma^3_{44}=-\cos\theta \sin\theta.
 \end{cases}
\end{eqnarray}
\indent The non-zero components $R_{hijk}$ and $S_{ij}$ of $R$ and $S$ of the Bardeen metric \eqref{BM} along with its scalar curvature $k$ are :
\begin{eqnarray}\label{RRS}
\begin{cases}
R_{1212}=\frac{M(15\rho^2e^2-2\rho^4_1)}{\rho^7_1}, \  
R_{1313}=\frac{M\rho^2 (\rho^2_1-3e^2)(-2M\rho^2+\rho^3_1)}{\rho^8_1}=\frac{1}{\sin^2\theta}R_{1414}, \\

R_{2323}=-\frac{M\rho^2(\rho^2_1-3e^2)}{\rho^2_1(-2M\rho^2+\rho^3_1)}=\frac{1}{\sin^2\theta}R_{2424}, \\

R_{3434}=\frac{2M\rho^4 \sin^2\theta}{\rho^3_1};  \\
S_{11}=\frac{3Me^2(2\rho^2_1-3e^2)(-2M\rho^2+\rho^3_1)}{\rho^{10}_1},  \
S_{22}=\frac{3e^2M(5\rho^2-2\rho^2_1)}{\rho^4_1(-2M\rho^2+\rho^3_1)},  \\
S_{33}=-\frac{6Me^2\rho^2}{\rho^5_1}=\frac{1}{\sin^2\theta}S_{44}; \

 \text{and}  \
k=\frac{6Me^2(5\rho^2-4\rho^2_1)}{\rho^7_1}.
\end{cases}
\end{eqnarray}
\indent Let $L^1=(g\wedge g)$, $L^2=(g\wedge S)$ and $L^3=(S\wedge S)$. Then the non-zero components of these Kulkarni-Nomizu products are written as:
\begin{eqnarray}\label{KNP}
\begin{cases}
L^1_{1212}=2,\ L^1_{1313}=2(\rho^2-\frac{2M\rho^4}{\rho^3_1})=\frac{1}{\sin^2\theta}L^1_{1414}, \\
L^1_{2323}=\frac{2\rho^2 \rho^3_1}{2M\rho^2-\rho^3_1}=-\frac{1}{\sin^2\theta}L^1_{2424}, \
L^1_{3434}= -2\rho^2 \sin^2\theta;\\ \\
%
L^2_{1212}=\frac{6Me^2(5\rho^2-2\rho^2_1)}{\rho^7_1}, \
L^2_{1313}=-\frac{3Me^2\rho^2(4\rho^2_1-5\rho^2)(-2M\rho^2+\rho^3_1)}{\rho^{10}_1}=\frac{1}{\sin^2\theta}L^2_{1414}, \\
L^2_{2323}=\frac{3Me^2\rho^2(4\rho^2_1-5\rho^2)}{\rho^4_1(-2M\rho^2+\rho^3_1)}=\frac{1}{\sin^2\theta}L^2_{2424}, \
L^2_{3434}=\frac{12Me^2\rho^4\sin^2\theta}{\rho^5_1}; \\ \\
%
L^3_{1212}=\frac{18M^2e^4(2\rho^2_1-5\rho^2)^2}{\rho^{14}_1}, \

L^3_{1313}=\frac{36M^2e^4\rho^2(2\rho^2_1-5\rho^2)(-2M\rho^2+\rho^3_1)}{\rho^{15}_1}=\frac{1}{\sin^2\theta}L^3_{1414}, \\

L^3_{2323}=-\frac{36M^2e^4\rho^2(2\rho^2_1-5\rho^2)}{\rho^9_1(-2M\rho^2+\rho^3_1)}=\frac{1}{\sin^2\theta}L^3_{2424}, \

L^3_{3434}=-\frac{72M^2e^4\rho^4\sin^2\theta}{\rho^{10}_1}.
\end{cases}
\end{eqnarray} 
From \eqref{KNP} we can decompose Riemann tensor explicitly as follows:
\bea\label{RT}
R=\varrho_1 L^1 + \varrho_2 L^2 + \varrho_3 L^3
\eea
where $\varrho_1=\frac{M(18\rho^2_1-25\rho^2)}{25\rho^2 \rho^3_1}$, $\varrho_2=\frac{\rho^2_1(6\rho^2_1-5\rho^2)}{25e^2\rho^2}$ and $\varrho_3=\frac{(3\rho^2_1-5\rho^2)\rho^7_1}{150Me^4\rho^2}$. Contracting the relation \eqref{RT} the following relation is entailed:
\bea\label{Ein2}
S^2+ \beta S + \bar{\beta} g=0
\eea
where $\beta=\frac{3Me^2(4\rho^2_1-5\rho^2)}{\rho^7_1}$ and $\bar{\beta}=\frac{18M^2e^4(2\rho^2_1-5\rho^2)}{\rho^{12}_1}$.
\begin{pr}\label{pr1}
	The Bardeen metric \eqref{BM} is not Einstein manifold but it is (i) Einstein manifold of level $2$ as well as (ii) fulfills Roter type condition.
\end{pr}
\begin{cor}\label{cr1}
	Since the Bardeen metric is Roter type, from Theorem $6.7$ of \cite{DGHS11} we obtain the following geometric structures of it:
	\begin{enumerate}[label=(\roman*)]
		\item $R\cdot R= \varrho_R Q(g,R)$,\ \ $\varrho_R=\frac{1}{2\varrho_2^3}(2((\varrho_2)^2-4\varrho_3\varrho_1)-2\varrho_1)=-\frac{M(2\rho^2_1-3\rho^2)^{5/2}}{\rho^2_1},$
		\item $R\cdot C=\varrho_R Q(g,C)$,
		\item $C\cdot R=\varrho_C Q(g,R)$, \ \ $\varrho_C=\varrho_R - (\frac{\kappa}{3}+\frac{\varrho_2}{2\varrho_3})+\frac{1}{4\varrho_3}=-\frac{M\rho^2(3\rho^2_1-5\rho^2)}{2\rho^7_1}$,
		\item $C\cdot C=\varrho_CQ(g,C)$,
		\item $R\cdot R=Q(S,R)+\varrho Q(g,C)$, \ \ $\varrho=\varrho_R+\frac{\varrho_2}{2\varrho_3}=\frac{2M(6\rho^2_1-7\rho^2)}{(3\rho^2_1-5\rho^2)\rho^3_1}$.
	\end{enumerate}
\end{cor}
The conformal curvature components of the metric \eqref{BM} are:
\begin{eqnarray}\label{CC}
\begin{cases}
C_{1212}=\frac{M\rho^2(3\rho^2_1-5\rho^2)}{\rho^7_1}=-\frac{1}{\rho^4 \sin^2\theta} C_{3434}; \\
C_{1313}=-\frac{M\rho^4(3\rho^2_1-5\rho^2)(-2M\rho^2+\rho^3_1)}{2\rho^5_1}=\frac{1}{ \sin^2\theta} C_{1414}; \\

C_{2323}= \frac{M\rho^4(3\rho^2_1-5\rho^2)}{2\rho^4_1(-2M\rho^2+\rho^3_1)}=\frac{1}{ \sin^2\theta} C_{2424}; \\

\end{cases}
\end{eqnarray}
Let $N^1_{abcd,f}=\nabla_fR_{abcd}$ and $N^2_{abcd,f}=\nabla_fC_{abcd}$. Then the non-vanishing components of the covariant derivatives of $R$ and $C$ are calculated and presented as below:
\begin{eqnarray}\label{CovR}
\begin{cases}
N^1_{1212,2}=\frac{3M\rho(12e^4-21e^2\rho^2+2\rho^5)}{\rho^9_1}; \\
N^1_{1213,3}=-\frac{3M\rho^3 (5\rho^2-4\rho^2_1)(-2M\rho^2+\rho^3_1)}{\rho^5_1}=-N^1_{1313,2}=\frac{1}{\sin^2\theta}N^1_{1214,4}=\frac{1}{\sin^2\theta}N^1_{1414,2}; \\

N^1_{2323,2}=\frac{3M\rho^3(5\rho^2-4\rho^2_1)}{\rho^4_1(-2M\rho^2+\rho^3_1)}=\frac{1}{\sin^2\theta}N^1_{2424,2}; \\

N^1_{2334,4}=\frac{3M\rho^5\sin^2\theta}{\rho^5_1}= -N^1_{2434,3}=-2N^1_{3434,2}; \\

\end{cases}
\end{eqnarray}
\begin{eqnarray}\label{CovC}
\begin{cases}
N^2_{1212,2}= \frac{\rho M(6e^4-23e^2\rho^2+6\rho^5)}{\rho^9_1}=-\frac{1}{\rho^4 \sin^2\theta} N^2_{3434,2}; \\
N^2_{1213,3}= \frac{3M\rho^3(3\rho^2_1-5\rho^2)(-2M\rho^2+\rho^3_1)}{2\rho^5_1}=\frac{1}{ \sin^2\theta} N^2_{1214,4}; \\

N^2_{1313,2}= -\frac{M\rho^3(6e^4-23e^2\rho^2+6\rho^4)(-2M\rho^2+\rho^3_1)}{2\rho^6_1}=-\frac{1}{ \sin^2\theta} N^2_{1414,2}; \\

N^2_{2323,2}=\frac{M\rho^3(6e^4-23e^2\rho^2+6\rho^4)}{2\rho^6_1(-2M\rho^2+\rho^3_1)}=\frac{1}{ \sin^2\theta} N^2_{2424,2}; \\
N^2_{2334,4}=-\frac{3M\rho^5(3\rho^2_1-5\rho^2)\sin^2\theta}{2\rho^7_1}=-N^2_{2434,3}; \\

\end{cases}
\end{eqnarray}
From \eqref{KNP}, \eqref{CovR} and \eqref{CovC} we get the following:
\begin{pr}\label{pr2}
	The Bardeen metric \eqref{BM} admits the following geometric structures:
	\begin{enumerate}[label=(\roman*)]
	\item $\nabla R=\Pi\otimes R + A\otimes (S\wedge S)$ where $\Pi=\left\lbrace 0, \frac{6\rho(8M-5\rho_1)}{5(-2M\rho^2+\rho^3_1)},0,0\right\rbrace $ and \\  $A=\left\lbrace 0, -\frac{\rho(29e^4+e^2(53\rho^2-8M\rho_1)+24(\rho^2-2M\rho_1 \rho^2))}{30M(-2M\rho^2+\rho^3_1)},0, 0 \right\rbrace $	,
	\item $\nabla R=A\otimes (g\wedge S)$ where $A=\left\lbrace 0,\frac{2\rho(8M-5 \rho_1)}{5(-2M\rho^2+\rho^3_1)},0,0  \right\rbrace $,\\
	\item \[\mathop{\mathcal{S}}_{\zeta_1,\zeta_2,\zeta_3}\nabla_{\zeta_1}C(\zeta_2,\zeta_3,X,Y)=\mathop{\mathcal{S}}_{\zeta_1,\zeta_2,\zeta_3}A(\zeta_1)\otimes C(\zeta_2,\zeta_3,X,Y)\] $$\text{where} \ \ A=\left\lbrace 0, \frac{5e^2(3\rho^2_1-7\rho^2)}{\rho \rho^2_1(3\rho^2_1-5\rho^2)},0,0 \right\rbrace .$$
	\end{enumerate}
\end{pr}
Let $B^1=R\cdot C$, $B^2=C\cdot R$, $F^1=Q(g,R)$, $F^2=Q(S,R)$, $F^3=Q(g,C)$ and $F^4=Q(S,C)$. Then, the non-vanishing components of the tensors $B^1$, $B^2$, $F^1$, $F^2$, $F^3$ and $F^4$ are given by (upto symmetry):
\begin{eqnarray}\label{RC}
\begin{cases}
B^1_{122313}=-\frac{3M^2\rho^4}{2\rho^{12}_1}(2\rho^2_1-3\rho^2)(3\rho^2_1-5\rho^2)=\frac{1}{\sin^2\theta}B^1_{122414}=-B^1_{121323}=-\frac{1}{\sin^2\theta}B^1_{121424},\\
B^1_{143413}=-\frac{3M^2\rho^6}{2\rho^{15}_1}(2\rho^2_1-3\rho^2)(3\rho^2_1-5\rho^2)(-2M\rho^2+\rho^2_1)\sin^2\theta=-B^1_{133414},\\
B^1_{243423}=\frac{3M^2\rho^6(2\rho^2_1-3\rho^2)(3\rho^2_1-5\rho^2)\sin^2\theta}{2\rho^9_1(-2M\rho^2+\rho^3_1)}=-B^1_{233424};
\end{cases}
\end{eqnarray}
\begin{eqnarray}\label{CR}
\begin{cases}
B^2_{1223,13}=  \frac{3M^2\rho^6(-4\rho^2_1+5\rho^2)(3\rho^2_1-5\rho^2)}{2\rho^{14}_1}=-B^2_{1213,23}=\frac{1}{\sin^2\theta}B^2_{1224,14}=-\frac{1}{\sin^2\theta}B^2_{1214,24}, \\

B^2_{1434,13}= \frac{3M^2\rho^8}{2\rho^{15}_1}(3\rho^2_1-5\rho^2)(-2Mr^2+\rho^3_1)=-B^2_{1334,14}, \\

B^2_{2434,23}=-\frac{3M^2\rho^8(3\rho^2_1-5\rho^2)}{2\rho^9_1(-2M\rho^2+\rho^3_1)}=-B^2_{2334,24};	
\end{cases}
\end{eqnarray}
\begin{eqnarray}\label{qgR}
\begin{cases}
F^1_{1223,13}= -\frac{3M\rho^4}{\rho^7_1}(-4\rho^2_1+5\rho^2)=-\frac{1}{\sin^2\theta}F^1_{1214,24}=\frac{1}{\sin^2\theta}F^1_{1224,14}=-F^1_{1213,23}, \\
F^1_{1434,13}= \frac{3M\rho^6}{\rho^4_1}(2M\rho^2-\rho^3_1)\sin^2\theta=-F^1_{1334,14}, \\

F^1_{2434,23}= \frac{3M\rho^6 \sin^2\theta}{\rho^2_1(-2M\rho^2+\rho^3_1)}=-F^1_{2334,24};
\end{cases}
\end{eqnarray}
\begin{eqnarray}\label{qSR}
\begin{cases}
F^2_{1223,13}=\frac{3M^2e^2\rho^4}{\rho^{12}_1}(14\rho^2_1+13\rho^2)=-F^2_{1213,23}=\frac{1}{\sin^2\theta}F^2_{1224,14}=-\frac{1}{\sin^2\theta}F^2_{1214,24}, \\

F^2_{1434,13}= -\frac{12M^2e^2\rho^6}{\rho^{13}_1}(-2M\rho^2+\rho^3_1)\sin^2\theta=-F^2_{1334,14}, \\

F^2_{2434,23}= \frac{12M^2e^2\rho^6\sin^2\theta}{\rho^7_1(-2M\rho^2+\rho^3_1)}=-F^2_{2334,24};
\end{cases}
\end{eqnarray}
\begin{eqnarray}\label{qgC}
\begin{cases}
F^3_{1223,13}=\frac{3M\rho^4}{2\rho^7_1}(3\rho^2_1-5\rho^2)=-F^3_{1213,23}=\frac{1}{\sin^2\theta}F^3_{1224,14}=-\frac{1}{\sin^2\theta}F^3_{1214,24}, \\

F^3_{1434,13}= \frac{3M\rho^6}{2\rho^{10}_1}(3\rho^2_1-5\rho^2)(-2M\rho^2+\rho^3_1)\sin^2\theta=-F^3_{1334,14}, \\

F^3_{2434,23}=- \frac{3M\rho^6(3\rho^2_1-5\rho^2)\sin^2\theta}{2\rho^4_1(-2M\rho^2+\rho^3_1)}=-F^3_{2334,24};
\end{cases}
\end{eqnarray}
\begin{eqnarray}\label{qSC}
\begin{cases}
F^4_{1223,13}=\frac{3M^2e^2\rho^4}{2\rho^{14}_1}(3\rho^2_1-5\rho^2)(6\rho^2_1-7\rho^2)=\frac{1}{\sin^2\theta} F^4_{1224,14}=-F^4_{1213,23}=\frac{1}{\sin^2\theta} F^4_{1214,24}, \\
F^4_{1434,13}=-\frac{3M^2e^2\rho^6}{\rho^{17}_1}(3\rho^2_1-5\rho^2)^2(-2M\rho^2+\rho^3_1)\sin^2\theta=-F^4_{1334,14}, \\

F^4_{2434,23}=\frac{3M^2e^2\rho^6(3\rho^2_1-5\rho^2)}{\rho^{11}_1(-2M\rho^2+\rho^3_1)}=F^4_{2334,24}.
\end{cases}
\end{eqnarray}
\begin{pr}\label{pr3}
From \eqref{RC}-- \eqref{qSC} we obtain the following pseudosymmetric type curvature relations for the metric \eqref{BM}:
\bea
C\cdot R -R\cdot C=\bar{\varrho}_2\ Q(S,R) + \bar{\varrho}_1\ Q(g,R)
\eea
where $\bar{\varrho}_1=-\frac{M(3\rho^2_1-5\rho^2)(\rho^2(6\rho^2_1-7\rho^2)-(2\rho^2_1-3\rho^2)^2)}{2(6\rho^2_1-7\rho^2)\rho^7_1}$ \ \ $\bar{\varrho}_2=1-\frac{3}{14}e^2(\frac{12}{6\rho^2_1-7\rho^2}+\frac{5}{\rho^2_1})$
\bea
\text{and} \ \ \ \ \ C\cdot R -R\cdot C=Q(S,C) + \bar{\varrho}_3\ Q(g,C) 
\eea
where $\bar{\varrho}_3=\frac{2M(4\rho^2_1-5\rho^2)e^2}{\rho^7_1}$.
\end{pr}

\indent Let $B^3=W\cdot R$ and $B^4=K\cdot R$. Then the non-zero components of the tensors $W \cdot R$ and $K\cdot R$ are given by :
\begin{eqnarray}\label{WR}
\begin{cases}
B^3_{1223,13}=\frac{3M^2\rho^6}{2\rho^{14}_1}(-4\rho^2_1+5\rho^2)(3\rho^2_1-5\rho^2)=-B^3_{1213,23}=\frac{1}{\sin^2\theta} B^3_{1224,14}=-\frac{1}{\sin^2\theta} B^3_{1214,24}, \\
B^3_{1434,13}=\frac{3M^2\rho^8}{2\rho^{15}_1}(3\rho^2_1-5\rho^2)(-2M\rho^2+\rho^3_1)\sin^2\theta=-B^3_{1334,14},\\

B^3_{2434,23}=-\frac{3M^2\rho^8(3\rho^2_1-5\rho^2)\sin^2\theta}{2\rho^9_1(-2M\rho^2+\rho^3_1)}=-B^3_{2334,24};
\end{cases}
\end{eqnarray}
\begin{eqnarray}\label{KR}
\begin{cases}
B^4_{1223,13}=-\frac{3M^2\rho^4}{2\rho^{14}_1}(-4\rho^2_1+5\rho^2)(8e^4-5e^2\rho^2+2\rho^4)=-B^4_{1213,23}=\frac{1}{\sin^2\theta} B^4_{1224,14}=-\frac{1}{\sin^2\theta} B^4_{1214,24}, \\
B^4_{1434,13}=-\frac{3M^2\rho^6}{2\rho^{15}_1}(8e^4-5e^2\rho^2+2\rho^4)(-2M\rho^2+\rho^3_1)\sin^2\theta=-B^4_{1334,14},\\

B^4_{2434,23}=\frac{3M^2\rho^6(8e^4-5e^2\rho^2+2\rho^4)\sin^2\theta}{2\rho^9_1(-2M\rho^2+\rho^3_1)}=-B^4_{2334,24};
\end{cases}
\end{eqnarray}
From \eqref{qgR}, \eqref{WR} and \eqref{KR} we get the following:
\begin{pr}\label{pr4}
	The Bardeen metric \eqref{BM} fulfills the curvature conditions
	$$W\cdot R=-\frac{M\rho^2(3\rho^2_1-5\rho^2)}{2\rho^7_1}Q(g,R)$$ and 
	$$K\cdot R=\frac{M(8e^4-5e^2\rho^2+2\rho^4)}{2\rho^7_1}Q(g,R).$$
	
\end{pr}
From propositions \ref{pr1}--\ref{pr4} and  from corollary \ref{cr1} we can conclude about the curvature properties of the Bardeen spacetime metric \eqref{BM} as follows:
\begin{thm}
	The Bardeen metric \eqref{BM} admits the following curvature restricted geometric structures:
	\begin{enumerate}[label=(\roman*)]
	\item It is a pseudosymmetric spacetime and as a result it realizes Ricci pseudosymmetry, conformal pseudosymmetry, concircular pseudosymmetry, conharmonic pseudosymmetry and projective pseudosymmetry,
	\item it is also pseudosymmetric due to conformal curvature, concircular curvature and conharmonic curvature,
	\item the difference tensor $C\cdot R - R\cdot C$ is linearly dependent with the tensors $Q(g,R)$ and $Q(S,R)$ as well as it is also linearly dependent with the tensors $Q(g,C)$ and $Q(S,C)$,
	\item the pseudosymmetric type condition $$R\cdot R- Q(S,R)=\varrho\ Q(S,C)$$ is possessed by this spacetime and it is also equipped with pseudosymmetric Weyl tensor,
	where $\varrho$ is a smooth scalar function given in corollary \ref{cr1},
	\item it is a weakly generalized recurrent manifold satisfying special recurrent like structure $\nabla R=A\otimes (g\wedge S)$,
	\item its conformal $2$ forms are recurrent,
	\item it is a Roter type spacetime and is an Einstein manifold of level $2$,
	\item Ricci tensor is Weyl compatible as well as Riemann compatible.
	\end{enumerate}
	
\end{thm}

\begin{rem}
	From the components of various curvatures we conclude that the following geometric structures are not admitted by the Bardeen metric \eqref{BM}:
		\begin{enumerate}[label=(\roman*)]
	\item any semisymmetric type conditions for $C, P, W, K, S$,
	\item Ricci generalized pseudosymmetry,  
	\item $D$-Venzi space for $E=C, R, P,W, K,$ 
	\item Einstein or quasi-Einstein condition,
	\item curvature $2$-forms recurrence, 
	\item Codazzi type Ricci tensor or cyclic parallel Ricci tensor,
	\item Chaki pseudosymmetry, 
	\item Weak symmetry.
\end{enumerate}  
\end{rem}

\section{\bf Energy momentum tensor of Bardeen black hole metric}

\indent In general theory of relativity, Einstein describes the physics of a spacetime in terms of geometry  by the system of equations 
$$S-\frac{k}{2}g+\Lambda g= \frac{8\pi G}{c^4}T$$
where $k$ is the scalar curvature, $S$ the Ricci curvature and $T$ the Energy momentum tensor of the spacetime. Also $\Lambda$ represents the cosmological constant, $G$ is the gravitational constant, $c$ is the speed of light in vacuum.

\indent Taking $\frac{8\pi G}{c^4}=1$ the components of Energy momentum tensor are given by: 
$$\begin{array}{c}
 T_{11}=-\frac{(-2M\rho^2+\rho^3_1)(6Me^2+\rho^5_1\Lambda)}{\rho^8_1}, \\
 T_{22}=\frac{6Me^2+\rho^5_1\Lambda}{\rho^2_1 (-2M\rho^2+\rho^3_1)}, \\
 T_{33}=\frac{\rho^2(3Me^2(2\rho^2_1-5\rho^2)+\rho^7_1\Lambda)}{\rho^7_1}=\frac{1}{\sin^2\theta}T_{44}.
\end{array}$$ 

\indent The non-zero components of the tensor $R \cdot T$ are
$$\begin{array}{c}
(R \cdot T)_{1313}=\frac{15M^2e^2\rho^4}{8\rho^{15}_1}(-2\rho^2_1+3\rho^2)(-2M\rho^2+\rho^3_1)=\frac{1}{sin^2\theta}(R \cdot T)_{1414}, \\

(R \cdot T)_{2323}=-\frac{15M^2e^2\rho^4(-2\rho^2_1+3\rho^2)}{8\rho^9_1(-2M\rho^2+\rho^3_1)}=\frac{1}{sin^2\theta}(R \cdot T)_{2424}, \\

\end{array}$$ 
\indent The non-zero components of the tensor $Q(g,T)$ are
$$\begin{array}{c}
 Q(g,T)_{1313}=\frac{15Me^2\rho^4}{8\rho^{10}_1}(-2M\rho^2+\rho^3_1)=\frac{1}{\sin^2\theta}Q(g,T)_{1414}, \\

 Q(g,T)_{2323}=-\frac{15Me^2\rho^4}{8\rho^4_1(-2M\rho^2+\rho^3_1)}=\frac{1}{\sin^2\theta}Q(g,T)_{2424}, \\
 
\end{array}$$ 
\indent Also the non-zero components of the tensor $C\cdot T$ are
$$\begin{array}{c}
(C\cdot T)_{1313}=-\frac{15Me^2\rho^6(3\rho^2_1-5\rho^2)(-2M\rho^2+\rho^3_1)}{2\rho^{17}_1}=\frac{1}{\sin^2\theta}(C\cdot T)_{1414}, \\
(C\cdot T)_{2323}=-\frac{15Me^2\rho^6(3\rho^2_1-5\rho^2)}{2\rho^{11}_1(-2M\rho^2+\rho^3_1)}=\frac{1}{\sin^2\theta}(C\cdot T)_{2424}, \\

\end{array}$$ 
 We can state the following, in view the above components:
\begin{thm}
	  The Bardeen spacetime \eqref{BM} has the energy momentum tensor concurring the following properties:
	\begin{enumerate}[label=(\roman*)]
	\item $R\cdot T=-\frac{M(2\rho^2_1-3\rho^2)}{\rho^5_1}Q(g,T)$ i.e., the energy momentum tensor is pseudosymmetric,
	\item $C\cdot T=-\frac{M\rho^2(3\rho^2_1-5\rho^2)}{2\rho^7_1}Q(g,T)$ i.e., the energy momentum tensor is pseudosymmetric due to Weyl tensor,
	\item the energy momentum tensor is Riemann compatible as well as Weyl compatible.	 
\end{enumerate}
\end{thm}

\section{\bf Bardeen black hole metric and Reissner-Nordstr\"om metric}
%
\indent Reissner-Nordstr\"om metric \cite{Kowa06} is a stationary solution of Einstein-Maxwell field equations with zero cosmological constant. Physically, it represents
the exterior gravitational field of a charged black hole. In spherical coordinates $(t,r,\theta, \phi)$, the Reissner-Nordstr\"om metric is given by \
\begin{eqnarray}\label{VBM}
ds^2=g_{ij}dx^idx^j=-\left(1-\frac{2m}{r}+\frac{q^2}{r^2}\right)dt^2-2dtdr+r^2(d\theta^2+\sin^2\theta d\phi^2),
\end{eqnarray}
where $t$ is the time coordinate, $r$ is the radial coordinate, the parameters $m$ is the mass of the body and $q$ is the charge of the body. We note that Bardeen spacetime is also a model of a charged black hole. Unlike Reissner-Nordstr\"om metric it does not bear curvature singularity. Hence we compare their curvature properties as this comparison compare a charge black hole with curvature singularity with curvature singularity free black hole, the Bardeen black hole.
\textbf{Similarities:}
\begin{enumerate}[label=(\roman*)]
	\item both the black holes are Roter type ,
	\item both spacetimes describe  $Ein(2)$ manifolds,
	\item both the black holes are pseudosymmetric,

	\item conformal curvature $2$-forms are recurrent for both,
	\item both the spacetimes have Riemann compatible and Weyl compatible Ricci tensor.
	
\end{enumerate} 
However, they have the following dissimilar properties:\\
\textbf{Dissimilarities:}
\begin{enumerate}[label=(\roman*)]
\item scalar curvature of  Reissner-Nordstr\"om spacetime vanishes while it doesn't for the Bardeen spacetime,
\item the Bardeen spacetime comes out with a weakly generalized recurrent manifold while  Reissner-Nordstr\"om spacetime doesn't,
\item also the Bardeen spacetime admits special recurrent like structure $\nabla R=A\otimes (g\wedge S)$ ($A$ being some $1$-form) but  Reissner-Nordstr\"om spacetime does not admit such recurrence.

\end{enumerate}

\section{Acknowledgment}
The third author greatly acknowledges to The University Grants Commission, Government of India for the award of Junior Research Fellow. All the algebraic computations of Section $3$ and $4$ are performed by a program in Wolfram Mathematica.
%


\end{document}